\documentclass[11pt,reqno]{amsart}

\usepackage{euscript}
\usepackage{pb-diagram}
\usepackage{lamsarrow}
\usepackage{pb-lams}
\usepackage{epsfig}
\usepackage{amssymb}

\theoremstyle{plain}

\theoremstyle{definition}

\DeclareMathOperator{\Hom}{{\rm Hom}}

\DeclareMathOperator{\rank}{{\rm rank}}

\newcommand{\ssw}{{\bf sw}}

\newcommand{\et}{\EuScript{T}}

\newcommand{\bL}{{\mathbb L}}

\newcommand{\cs}{\langle \chi_0\rangle}

\def\C{\mathbb C}
\def\Q{\mathbb Q}

\def\Z{\mathbb Z}

\def\im{{\rm Im}}

\newcommand{\x}{(X,0)}
\newcommand{\tx}{\tilde{X}}
\newcommand{\tz}{\tilde{Z}}
\newcommand{\call}{{\mathcal L}}

\newcommand{\calo}{{\mathcal O}}

\newcommand{\calt}{{\mathcal T}}

\newcommand{\calj}{{\mathcal J}}

\begin{document}

\title{Line bundles associated with normal surface singularities}	

\author{Andr\'as N\'emethi}
\address{Department of Mathematics\\Ohio State University\\Columbus, OH 43210}
\email{nemethi@math.ohio-state.edu}
\urladdr{http://www.math.ohio-state.edu/ \textasciitilde nemethi/}
\thanks{The author is partially supported by NSF grant DMS-0304759.}




\maketitle
\pagestyle{myheadings}
\markboth{{\normalsize A. N\'emethi}}{
{\normalsize Line bundles associated with normal surface singularities}}

{\small

\section{Introduction}

In \cite{SWI} L. Nicolaescu and the author
formulated a conjecture which relates the 
geometric genus of a complex analytic normal surface singularity 
$\x$ (whose link $M$ is a rational homology sphere)
 with the Seiberg-Witten invariant of $M$ associated with the ``canonical''
$spin^c$ structure of $M$. (The interested reader is invited to see 
the articles \cite{SWI,SWII,SWIII,NOSZ,NW,NWuj} for the verification of the 
conjecture in different cases.)
Since the Seiberg-Witten theory of the link $M$ provides a rational number
for {\em any} $spin^c$ structure (which are classified by $H_1(M,\Z)$),
it was a natural challenge to search for a complete set of conjecturally valid
identities, which involve all the Seiberg-Witten invariants (giving an 
analytic -- i.e. singularity theoretical --  interpretation of them).

The formulation of this set of identities is one of the goals of the 
present article. In fact (similarly as in \cite{SWI}), we formulate 
conjecturally valid inequalities  which became equalities in special 
rigid situations.
 In this way, the Seiberg-Witten invariants 
determine optimal topological upper bounds for the dimensions of the 
first sheaf-cohomology of line bundles living on some/any  resolution of $\x$.
Moreover, for $\Q$-Gorenstein singularities and some 
``natural'' line bundles equality holds. 

The first part of the article constructs these ``natural'' holomorphic
line bundles on the resolution $\tx$ of $\x$.  This construction
automatically provides a natural splitting of the exact sequence
$$0\to Pic^0(\tx)\to Pic(\tx)\stackrel{c_1}{\to} H^2(\tx,\Z)\to 0.$$
The line-bundle construction is compatible with abelian covers. This 
allows us to reformulate the conjecture in its second version
which relates the echivariant geometric genus of the universal 
(unbranched) abelian cover of $\x$ with the Seiberg-Witten invariants of 
the link.

In the last section we verify the conjecture for rational singularities.

\section{Preliminaries}

\subsection{}\label{2.1}  Let $\x$ be a complex analytic normal surface singularity.
We fix a good resolution $\pi : \tx \to X$ of $\x$ (over a small fixed Stein
representative $X$ of $\x$) such that the exceptional divisor
$E:=\pi^{-1}(0)$ has only normal crossing singularities. We denote the irreducible 
components of $E$ by $\{E_j\}_{j\in \calj}$. 
The boundary $\partial \tx$ of $\tx$ can be identified with the link $M$ of $\x$,
which is an oriented connected 3-manifold. In this article, we will assume that $M$ is an
{\em rational homology sphere} (i.e. $H_1(M,\Q)=0$).  This implies that all the curves
$E_j$ are rational, and $H_1(\tx,\Z)=0$. 

The exact sequence of $\Z$-modules 
\begin{equation*}
0\to L \stackrel{i}{\to}  L'\to H\to 0
\tag{1}
\end{equation*}
will stay for the homological exact sequence 
\begin{equation*}
0\to H_2(\tx,\Z)\to H_2(\tx,M,\Z)\stackrel{\partial}{\longrightarrow}
 H_1(M,\Z)\to 0,
\tag{2}
\end{equation*}
or, via  Poincar\'e duality, for 
\begin{equation*}
0\to H^2_c(\tx,\Z)\to H^2(\tx,\Z)\to H^2(M,\Z)\to 0.
\tag{3}
\end{equation*}
$L$,  considered as in (2), 
 is freely generated by the homology classes $[E_j]_{j\in \calj}$.
For each $j$, consider a small transversal disc $D_j$ in $\tx$ with
$\partial D_j\subset \partial \tx$. Then $L'$ is freely generated by
the classes $[D_j]_{j\in \calj}$. We set $g_j:=\partial [D_j]=[\partial D_j]\in H$.
In the sequel, for simplicity we write $E_j,\ D_j$ for $[E_j]$ and $[D_j]$.

On $L$ there is a symmetric, non-degenerate, negative definite intersection 
form $(\cdot,\cdot )$ (with associated matrix $B_{ij}=(E_j,E_j)$), cf. \cite{Artin62,GRa}.
Then the morphism $i:L\to L'$  can be identified 
with $L\to \Hom(L,\Z)$ given by $l\mapsto (l,\cdot)$.
The intersection form has a natural extension to $L_\Q=L\otimes \Q$,
and it is convenient to identify $L'$ with a sub-lattice  of $L_\Q$:
$\alpha\in \Hom (L,\Z)$ corresponds with the unique $l_\alpha\in L_\Q$ which
satisfies $\alpha(l)=(l_\alpha,l)$ for any $l\in L$. By this identification, $D_j$,
considered in $L_\Q$ (and written in the base $\{E_j\}_j$), 
is the $j^{th}$ column  of $B^{-1}$, and $(D_j,E_i)=\delta_{ji}$. 

Elements $x=\sum_j r_jE_j \in L_\Q$ will be called (rational) cycles.  
If $x_i=\sum_jr_{j,i}E_j$ for $i=1,2$, then
$\min\{x_1,x_2\}:=\sum_j\min\{r_{j,1},r_{j,2}\}E_j$.
$x^2$ means  $(x,x)$; $e_j:=E_j^2$. 
We define the support $|x|$ of $x$ by $\cup E_j$, where the union runs over
$\{j:r_j\not=0\}$. 

\subsection{}\label{2.2} One can  consider two types of ``positivity conditions''
 for rational cycles.  The first one is considered in $L$:
A cycle $x=\sum_j  r_jE_j\in  L_\Q$ is  called {\em effective}, 
denoted by $x\geq 0$, if $r_j\geq 0$  for any $j$.  Their set is denoted by  
$L_{\Q,e}$, while $L_e' :=L_{\Q,e}\cap L'$ and
$L_e :=L_{\Q,e}\cap L$.  We write $x\geq  y$ if  
$x-y\geq 0$. $x>0$ means $x\geq 0$ but $x\not=0$. 
$\geq $ provides a partial order of $L_\Q$. 

The second is the {\em numerical effectiveness} of the rational cycles,
 i.e. positivity considered in $L'$. 
We define $NE_\Q:=\{x\in L_\Q : (x,E_j)\geq  0 \ \mbox{for any $j$}\}$.
In fact, $NE_\Q$ is the positive cone in $L_\Q$ generated by 
$\{D_j\}_j$, i.e. it is $\{\sum_j r_j D_j, r_j\geq 0 \ \mbox{for any $j$}\}$.
Since $B$ is negative definite, all the entries of $D_j$ are 
{\em strict}  negative.  In particular, $-NE_\Q\subset L_{\Q,e}$.

\subsection{Characteristic elements. $Spin^c$-structures.}\label{2.3} The 
set of characteristic elements are defined by 
$$Char=Char (\tx):=\{k\in L': \, (k,x)+(x,x)\in 2\Z \ \mbox{for any $x\in L$}\}.$$
There is a {\em canonical} characteristic element $K\in Char$
defined by $(K,E_j)=-(E_j,E_j)-2$  for any $j$. 

Clearly, $Char=K+2L'$. There is a natural
action of $L$ on $Char$ by $x*k:=k+2x$ whose orbits are of type
$k+2L$. Obviously, $H$ acts freely and transitively on the set of orbits  
by $[l']*(k+2L):=k+2l'+2L$ (in particular, they have the same cardinality). 

If $\tilde{X}$ is as above, then  the first Chern class (of the
associated determinant line bundle)
realizes an identification between the $spin^c$-structures
$Spin^c(\tilde{X})$ on $\tilde{X}$ and $Char\subset L'=
H^2(\tilde{X},\Z)$ (see e.g. \cite{GS}, 2.4.16). 
The restrictions to $M$ defines an identification 
of the $spin^c$-structures $Spin^c(M)$ of $M$ with the set of orbits of
$Char$ modulo $2L$; and this identification is compatible with the action of 
$H$ on both sets. In the sequel, we think about $Spin^c(M)$ by this 
identification, hence any $spin^c$-structure 
 of $M$ will be represented by an orbit  $[k]:=k+2L\subset Char$. 


\subsection{Liftings.}\label{2.4} If $H$ is not trivial,
then the exact sequence (1) does not split. 
Nevertheless, we will consider some 
``liftings'' (set theoretical sections) of the element of $H$ into $L'$. They correspond to
the positive cones in $L_\Q$ considered in \ref{2.2}.

More precisely, for any $l'+L=h\in H$, let $l'_e(h)\in L'$ be the unique minimal effective 
rational cycle in $L_{\Q,e}$ whose class is $h$. Clearly, the set 
$\{l'_e(h)\}_{h\in H}$ is exactly the closed/open  unit cube 
$Q:=\{\sum_j r_jE_j \in L'; 0\leq r_j<1\}$.

Similarly, for any $h=l'+L$, the intersection $(l'+L)\cap NE_\Q$
has a unique maximal element $l'_{ne}(h)$,
and the intersection $(l'+L)\cap (-NE_\Q)$)
has a unique minimal element $\bar{l}'_{ne}(h)$
 (cf. \cite{NOSZ}, 5.4).  By their definitions $\bar{l}'_{ne}(h)=-l'_{ne}(-h)$.

The elements $\bar{l}'_{ne}(h) $ were introduced in 
\cite{NOSZ} in order to compute the Ozsv\'ath-Szab\'o Floer homology of $M$
for some singularity links (their  notation in that paper 
is  $l'_{[k]}$, where  $[k]=K+2(l'+L)$). Using these elements, 
[loc. cit.] defines the distinguished representative $k_r$ of $[k]$ by 
$k_r:=K+2\bar{l}'_{ne}(h)$  ($h=l'+L$). 

For some $h$, $\bar{l}'_{ne}(h)$ 
might be situated in $Q$, but, in general, this is not the case
(cf. \ref{6.3} and \ref{6.4}). 
In general, the characterization of all the elements $\bar{l}'_{ne}(h)$ is  not simple
(for the cases when $M$ is a lens space or a Seifert manifold, see 
\cite{NOSZ}).

\subsection{RR}\label{RR} In some of our arguments, we use the same notation
for $l=\sum n_jE_j\in L $  and the algebraic 
cycle $\sum n_jE_j $ of $\tx$ supported by $E$.  
E.g., for any $l=\sum n_jE_j\in L $
one can take the line bundle $\calo_{\tx}(l):=\calo_{\tx}(\sum n_jE_j)$.
If $l>0$ set  $\calo_l:=\calo_{\tx}/\calo_{\tx}(-l)$, and let 
$\chi(l):=\chi(\calo_l)=h^0(\calo_l)-h^1(\calo_l)$ be its Euler-characteristic.
It can be computed combinatorially by Riemann-Roch: $\chi(l)=-(l,l+K)/2$. 

By analogy, for any $l'\in L_\Q$, we define $\chi(l'):=-(l',l'+K)/2$. 

\section{Line bundles on $\tx$.}

\subsection{}\label{3.1} Let $\pi:(\tx,E)\to \x$ be a fixed good resolution of
$\x$ as in section 2. Since $H^1(\tx,\Z)=0$, one has
\begin{equation*}
0\to Pic^0(\tx)\to Pic(\tx)\stackrel{c_1}{\to} L'\to 0,
\tag{1}
\end{equation*}
where $Pic^0(\tx)=H^1(\tx,\calo_{\tx})$, $Pic(\tx)= H^1(\tx, \calo^*_{\tx})$, and
$c_1({\call})=\sum_j \deg(\call|E_j)\, D_j$ is the set of Chern numbers 
(multidegree) of $\call$. 
Notice that $c_1(\calo_{\tx}(l))=l$, hence $c_1$ admits a group-section
$s_L:L\to Pic(\tx)$ above the subgroup $L$ of $L'$.  In the sequel we will construct 
a (natural) group section $s:L'\to Pic(\tx)$ of (1) which extends $s_L$
(and is compatible with abelian coverings). 

Clearly, if $Pic^0(\tx)=0$ (i.e. if $\x$ is rational, 
cf. \cite{Artin62})
then there is nothing to construct:
$c_1$ is an isomorphism, and $s:=c_1^{-1}$ identifies the lines bundles with
their 
multidegree. Nevertheless, for non-rational singularities, even the existence 
of any kind of splitting of (1) is not so obvious. 

\subsection{}\label{3.2} Notice that $\tx\setminus E\approx X\setminus \{0\}$ 
has the homotopy type of $M$, hence the abelianization map 
$\pi_1(\tx\setminus E)=\pi_1(M)\to H$ defines a regular Galois covering
of $\tx\setminus E$. This has a unique extension $c:Z\to \tx$ with $Z$ normal
and $c$ finite \cite{GR}. The (reduced) branch locus of $c$ is included in $E$, 
and the Galois action of $H$ extends to $Z$ as well. Since $E$ is a  normal
crossing divisor, the only singularities what $Z$ might have are cyclic quotient
singularities (situated above $Sing(E)$). Let $r:\tz\to Z$ be  a resolution of 
these singular points of $Z$, such that $(c\circ r)^{-1}(E)$ is a normal crossing divisor. 

Notice that, in fact, $\tz$ is a good resolution of the universal 
(unbranched) abelian cover 
$(X_a,0)$ of $\x$. (Here, $(X_a,0)$ is the unique normal singular germ  corresponding
to the regular covering of $X\setminus \{0\}$ 
associated with $\pi_1(X\setminus \{0\})\to H$.) 
Hence, one can consider the following commutative diagram:
$$\begin{array}{ccccccccc}
0&\to&L& \to & L'& \to & H & \to & 0\\
 & & \Big\downarrow & & 
\Big\downarrow \vcenter{%
\rlap{$\scriptstyle{p'}$}} 
& &
\Big\downarrow \vcenter{%
\rlap{$\scriptstyle{p_H}$}} 
& & \\
0&\to&L_a& \to & L_a'& \to & H_a & \to & 0 
\end{array} $$
Here, the first, resp. second,  horizontal line is  the exact sequence \ref{2.1}(3) 
applied for the resolution $\tx\to \x$, resp. for $\tz\to (X_a,0)$. The vertical arrows 
(pull-back of cohomology classes) are induced by $p=c\circ r$. 

\subsection{Lemma.}\label{3.3} $p_H=0$. {\em In particular, 
$p'(L')\subset L_a$ (i.e. any element $p'(l')$, $l'\in L'$, 
can be represented by a divisor 
supported by the exceptional divisor in $\tz$).}
\begin{proof} Denote by $M_a$ the link of $(X_a,0)$. The morphism
$p_H:H^2(M,\Z)\to H^2(M_a,\Z)$ is dual to $p_*:H_1(M_a,\Z)\to 
H_1(M,\Z)$, which is zero since $H_1(M_a,\Z)$ is the abelianization of 
the commutator subgroup of $\pi_1(M)$.\end{proof}

\noindent This shows that for any $l'\in L'$, one can take 
$\calo_{\tz}(p'(l'))\in Pic(\tz)$. 

\subsection{Theorem.}\label{prop3} {\em 
The  line bundle $\calo_{\tz}(p'(l'))$
is a pull-back of a unique element of $Pic(\tx)$.}

\begin{proof} We brake the proof into several steps.
Let $f:S\to T$ be one of the maps $r:\tz\to Z$, $c:Z\to \tx$ or
$p:\tz\to \tx$. In each case one has a commutative diagram of type

$$\begin{array}{ccccccccc}
H^1(T,\Z) & \to & Pic^0(T) & \to &  Pic(T) & \to & H^2(T,\Z) & \to & 0 \\
\Big\downarrow \vcenter{%
\rlap{$\scriptstyle{p''}$}} 
 & & \Big\downarrow   \vcenter{%
\rlap{$\scriptstyle{p^{0,*}}$}} 
& & 
\Big\downarrow \vcenter{%
\rlap{$\scriptstyle{p^*}$}} 
& &
\Big\downarrow \vcenter{%
\rlap{$\scriptstyle{p'}$}} 
& & \\
 H^1(S,\Z) & \to & Pic^0(S) & \to &  Pic(S) & \to & H^2(S,\Z) & \to & 0 
\end{array}$$
(I). {\em Assume that $f=c$.} Then:

$\bullet$ \ {\em $c'$ is injective } since 
$c'\otimes \Q: H^2(\tx,\Q)\approx H^2(Z,\Q)^H\hookrightarrow H^2(Z,\Q)$ is so
and $H^2(\tx,\Z)$ is free. 

$\bullet$ \ {\em $c^{0,*}$ is injective with image 
$Pic^0(Z)^H$.} Indeed (cf. also with the proof of \ref{3.9}), $c_*\calo_Z$
has a direct sum decomposition $\oplus_\chi \call_\chi$ into line bundles 
$\call_\chi$, where the sum
runs over all the characters of $H$, and for the 
trivial character $\call_1=\calo_{\tx}$.
Therefore, since $c$ is finite, $H^1(Z,\calo_Z)=H^1(\tx, c_*\calo_Z)=
\oplus_\chi H^1(\tx, \call_\chi)$, whose $H$-invariant part 
is $H^1(\tx, \calo_{\tx})$. 

$\bullet$ \ $\im(H^1(Z,\Z)\to Pic^0(Z))\cap \im p^{0,*}=0$. Indeed,
since any element of $\im p^{0,*}$ is  $H$-invariant, 
the above intersection  is in 
$\im(H^1(Z,\Z)^H\to Pic^0(Z)^H)$. But $H^1(Z,\Z)^H$ embeds into
$H^1(Z,\Q)^H=H^1(\tx,\Q)=0$, hence it is trivial. 

$\bullet$ {\em $p^*$ is injective,} a fact which follows from the 
previous three statements. 

\noindent (II). {\em Assume that $f=r$.}  Then $r''$ is an isomorphism
(since such a resolution does not modify $H_1(\cdot ,\Z)$ of the 
exceptional divisors),
$r^{0,*}$ is an isomorphism (since a quotient singularity has 
geometric genus zero), and $r'$ is injective. Hence (by a diagram 
check) $p^*$  is also injective.

\noindent (III). {\em Assume that $f=p=c\circ r$.} Then: 

$\bullet$ \ {\em For any $l'\in L'$, $\calo_{\tz}(p'(l'))$ is in the image of $p^*$.}
Indeed, take a line bundle $\call\in Pic(\tx)$ with $c_1\call=l'$.
Then $\call':=\calo_{\tz}(p'(l'))\otimes p^*\call^{-1}$
has trivial multidegree, and  it is in $Pic^0(\tz)^H$. 
But using (I-II)  $p^{0,*}$ is 
onto on $Pic^0(\tz)^H$, hence $\call'$ is a pull-back. Hence 
$\calo_{\tz}(p'(l'))$ itself is a pull-back.

$\bullet$ \ Again, from (I-II), {\em $p^*$ is injective}; a fact 
which ends the proof of \ref{prop3}. \end{proof}

\subsection{Notation.}\label{3.7} Write  $\calo_{\tx}(l')$
for the unique line bundle $\call\in Pic(\tx)$
with $p^*(\call)=\calo_{\tz}(p'(l'))$.

\vspace{2mm}

\noindent The proof of \ref{prop3} also implies the following fact.

\subsection{Corollary.}\label{3.8} {\em $s:L'\to Pic(\tx)$ defined by
$l'\mapsto \calo_{\tx}(l')$  is a group section of the exact sequence (1)
which  extends $s_L$ (cf. \ref{3.1}). }


\vspace{2mm}

The following result will illuminate a different aspect of the line 
bundles $\calo_{\tx}(l')$. In fact, the next proposition could also serve
as the starting point of a different construction of these line bundles. 

Below we will write $\hat{H}$ for the Pontrjagin dual $\Hom(H,S^1)$ 
of $H$. Recall that the natural map 
$$\theta:H\to \hat{H}, \ \ \mbox{induced by} \ \   [l']\mapsto e^{2\pi i (l',\cdot )}$$
is an isomorphism.

\subsection{Theorem.}\label{3.9} {\em Consider the finite covering
$c:Z\to \tx$, and set
$Q:=\{\sum r_jE_j\in L':\, 0\leq r_j<1\}\subset L'$ (cf. \ref{2.4})
 as above.  Then the $H$-eigenspace decomposition
of $c_*\calo_Z$ has the form: 
$$c_*\calo_Z=\oplus _{\chi\in \hat{H}}\call_\chi,$$
where $\call_{\theta(h)}=\calo_{\tx}(-l'_e(h))$  for any $h\in H$. In particular,}
$c_*\calo_Z=\oplus_{l'\in Q}\calo_{\tx}(-l').$

\begin{proof} The proof is based on a similar statement of Koll\'ar 
valid for cyclic coverings, see e.g. \cite{Kollar}, \S 9.

First notice that $c_*\call_Z$ is free. Indeed, since all the singularities of $Z$ 
are cyclic quotient singularities, this fact follows from the corresponding
statement for cyclic Galois coverings, which was verified in \cite{Kollar}.
Moreover, above $\tx\setminus E$ the covering is regular (unbranched)
corresponding to the regular representation of $H$. Therefore, $\rank c_*\call_Z=
|H|$, and it has an eigenspace decomposition $\oplus _{\chi}\call_\chi$,
where all the characters $\chi\in \hat{H}$ appear, and $\call_\chi|\tx\setminus E$
is a line bundle.  By a similar reduction as above to the cyclic case, one gets that
$\call_\chi$ itself is a line bundle.  Moreover, $\call_1$ (corresponding to the 
trivial character 1) equals $\calo_{\tx}$. 

Next we identify $\call_\chi$ for any character.  Fix $\chi_0\in \hat{H}$. It generates
the subgruop $\cs$ in $\hat{H}$. Write $\hat{K}:=\hat{H}/\cs$.  The morphism 
$\pi_1(M)\to H\to \cs\hat{}$ \,  defines a well-defined  normal  space $Y$ as a
cyclic Galois $\cs\hat{}$-covering of $\tx$. Clearly, one has the natural maps
$Z\stackrel{f}{\longrightarrow} Y\stackrel{e}{\longrightarrow}\tx$ 
with $ e\circ f=c$. 
By similar argument as above 
$f_*\calo _Z=\oplus _{\xi\in \hat{K}}\bar{\call}_\xi$
for some $\bar{\call}_\xi\in Pic(Y)$ and $\bar{\call_1}=\calo_Y$.
Since $e_*\bar{\call}_\xi=\oplus _{[\chi]=\xi} \call_\chi,$
and $[\chi_0]=0$ in $\hat{K}$, one gets that $\call_{\chi_0}$ is one of the 
summands of $e_*\bar{\call}_1=e_*\calo_Y$. In particular, $\call_{\chi_0}$
can be recovered from the cyclic covering $e$ as the $\chi_0$-eigenspace
of $e_*\calo_Y$. 

Assume that the order of $\chi_0$ is $n$, i.e. $\cs=\Z_n$, and we regard
$\chi_0$ as the distinguished generator of $\Z_n$.
Then for any $j\in \calj$,
$\chi_0(g_j)$ has the form $e^{2\pi im_j/n}$ for some (unique) 
$0\leq m_j<n$ (for the definition of $g_j$ see \ref{2.1}).  Using these integers, define
the divisor $B:=\sum_jm_jE_j$. Since for any fixed $i\in \calj$ one has
$$1=\chi_0([E_i])=\chi_0(\textstyle{\sum}_jB_{ij}g_j)=e^{2\pi i(\sum_jm_jE_j,E_i)/n},$$
one gets that $B/n\in L'$. Let $h=[B/n]$ be its class in $H$. Clearly,
$\theta(h)=\chi_0$ since
$$\theta(h)(g_j)=e^{2\pi i(\sum m_jE_j/n,D_j)}=e^{2\pi im_j/n}=\chi_0(g_j).$$
On the other hand, by \cite{Kollar}, 9.8 (and from the fact that all
the coefficients of $B/n$ are in the interval $[0,1)$), the $\chi_0$-eigenspace of 
$e_*\calo_Y$  is some line bundle $\call^{-1}$ with the properties 
(i) $\call^{\otimes n}=\calo_{\tx}(B)$ and 
(ii) $e^*\call=\calo_Y(e'(B/n))$. From (i) follows that  $c_1(\call)=B/n$, 
hence (ii), via our definition \ref{3.7}, reads as 
$\call=\calo_{\tx}(B/n)$. Therefore, $\call^{-1}=\calo_{\tx}(-B/n)$. 
Since $B/n$ is in the unit cube, $B/n=l'_e(h)$. 
\end{proof}

\section{Some cohomological computations}

\subsection{}\label{4.1} Let $\x$ and $\pi:\tx \to X$ as above. In this section
we analyse $h^1(\tx,\call)$ for any $\call\in Pic(\tx)$. 

The next \ref{4.3} is an improvement of the following
general (Kodaira, or Grauert-Riemenschneider type) vanishing
theorem (cf. \cite{Reid}, page 119, Ex. 15):
{\em If $c_1(\call)\in K+NE_\Q$, 
then $h^1(l,\call|_l)=0$ for any $l\in L$, 
$l>0$, hence $h^1(\tx,\call)=0$. }

\subsubsection{}\label{4.3} {\em Assume that $\x$ is a rational singularity.
If $c_1(\call)\in NE_\Q$,
then $h^1(l,\call|_l)=0$ for any $l>0$, 
$l\in L$, hence $h^1(\tx,\call)=0$ too. }

\vspace{2mm}

For the convenience of the reader we sketch a  proof. 
For any $l>0$ there exists $E_j\subset |l|$ such that $(E_j,l+K)<0$.
Indeed, $(E_j,l+K)\geq 0$ for any $j$ would imply $\chi(l)=-(l,l+K)/2\leq 0$, which 
would contradict the rationality of $\x$ \cite{Artin62}. Then,
 from the cohomology exact sequence of
$$0\to \call\otimes \calo_{E_j}(-l+E_j)\to \call|_l\to \call|_{l-E_j} \to 0$$
one gets $h^1(\call|_l)=h^1(\call|_{l-E_j})$, hence by induction $h^1(\call|_l)=0$.
We will generalize this  proof as follows.

\subsection{Proposition.}\label{4.4} {\em Let $\tx\to X$ be a good resolution
of a normal singularity $\x$ as above.

(a) For any $l'\in L'$ there exists a unique minimal element $l_{l'}\in L_e$ with 
$e(l'):=l'-l_{l'}\in NE_\Q$. 

(b) $l_{l'}$ can be found by the following (generalized Laufer's)  algorithm.
One constructs a ``computation sequence'' $x_0, x_1, \ldots, x_t\in L_e$
with $x_0=0$ and $x_{i+1}=x_i+E_{j(i)}$, where the index $j(i)$ is determined by the 
following principle. Assume that $x_i$ is already constructed.
Then, if $l'-x_i\in NE_\Q$, then one stops, and $t=i$. Otherwise, there exists 
at least one $j$ with $(l'-x_i,E_j)<0$. Take for $j(i)$ one of these $j$'s. 
Then this algorithm stops after a finitely many steps, and $x_t=l_{l'}$. 

(c) For any $\call\in Pic(\tx)$ with $c_1(\call)=l'$ one has:
$$h^1(\call)=h^1(\call\otimes \calo_{\tx}(-l_{l'}))-(l',l_{l'})-\chi(l_{l'}).$$
In particular (since $c_1(\call\otimes \calo_{\tx}(-l_{l'}))\in NE_\Q$),
the computation of any $h^1(\call)$ can be reduced (modulo the combinatorics of
$(L,(\cdot,\cdot))$) to the computation of  some $h^1(\call')$ with $c_1(\call')\in NE_\Q$. }

\vspace{2mm}

\noindent [Although the above decomposition $l'=e(l')+l_{l'}$ with
$l_{l'}\geq 0$ 
and $e(l')\in NE_\Q$ looks  similar to the Zariski decomposition
(cf. \cite{Sakai}), in general it is  different from it.]
\begin{proof} (a) First notice that since $B$ is negative definite, 
there exists at least one effective cycle $l$  with $l'-l\in NE_\Q$
(take e.g. a large multiple of some $x$ with $(x,E_j)<0$ for any $j$).
Next, we  prove that if $l'-l_i\in NE_\Q$ for $l_i\in L_e$,
$i=1,2$, and $l:=\min\{l_1,l_2\}$, then $l'-l\in NE_\Q$ as well.
For this, write $x_i:=l_i-l\in L_e$. Then $|x_1|\cap|x_2|=\emptyset$, 
hence for any fixed $j$, $E_j\not\subset |x_i|$ for at least one of the $i$'s.
Therefore, $(l'-l,E_j)=(l'-l_i,E_j)+(x_i,E_j)\geq 0$.

(b) First we prove that $x_i\leq l_{l'}$ for any $i$. For $i=0$ this is clear.
Assume that it is true for some $i$ but not for $i+1$, i.e. $E_{j(i)}\not\subset
|l_{l'}-x_i|$. But this would imply $(l'-x_i,E_{j(i)})=(l'-l_{l'},E_{j(i)})+(l_{l'}-x_i,E_{j(i)})\geq0$,
a contradiction. 
The fact that $x_i\leq l_{l'}$ for any $i$ implies that the algorithm must stop, and 
$x_t\leq l_{l'}$. But then by the minimality of $l_{l'}$ (part a) $x_t=l_{l'}$. (Cf. \cite{Laufer72}.)

(c) For any $0\leq i<t$, consider the exact sequence 
$$0\to \call\otimes \calo_{\tx}(-x_{i+1}) \to
\call\otimes \calo_{\tx}(-x_{i})
\to 
\call\otimes \calo_{E_{j(i)}}(-x_{i})
\to 0.$$
Since $\deg(\call\otimes \calo_{E_{j(i)}}(-x_{i}))=(l'-x_i,E_{j(i)})<0$, one gets
$h^0(\call\otimes \calo_{E_{j(i)}}(-x_{i}))=0$. 
Therefore
$$h^1(\call\otimes \calo_{\tx}(-x_{i}))-
h^1(\call\otimes \calo_{\tx}(-x_{i+1}))=
-\chi(\call\otimes \calo_{E_{j(i)}}(-x_{i}))$$
which equals $-(l',x_{i+1}-x_i)+\chi(x_i)-\chi(x_{i+1})$. 
Hence the result follows by induction.
\end{proof}

\subsection{Examples.}\label{4.5}  If $\call=\calo_{\tx}(l')$ for some $l'\in L'$
(cf. \ref{3.7}) then \ref{4.4}(c) reads as 
$$h^1(\calo_{\tx}(l')))=h^1(\calo_{\tx}(e(l')))-(l',l_{l'})-\chi(l_{l'}).$$
Additionally, if $\x$ is rational then 
$h^1(\calo_{\tx}(e(l')))=0$ by \ref{4.3}, hence
$h^1(\calo_{\tx}(l'))=-(l',l_{l'})-\chi(l_{l'})$. 
In particular, for $\x$ rational,
$h^1(\call)$ depends only on topological data and it is independent of
the analytic structure of $\x$.

\subsection{Definition.}\label{4.6} We will distinguish the following
set of rational cycles:
$$\bL':=\{ l'\in L':\, e(l')=l'_{ne}(h) \ \mbox{for some $h\in H$}\}=
\bigcup_{h\in H}l'_{ne}(h)+L_e.$$
It is easy to verify the following inclusions:
$L'_e\subset  \bigcup_{l'\in Q} -l'+L_e \subset \bL'$.

\section{The main conjecture}

\subsection{}\label{5.1}  In this section we present a conjecture which 
provides an optimal  upper bound for 
$h^1(\call)$ ($\call\in Pic(\tx)$) in terms of  the topology 
of $\x$ and $c_1(\call)$.
The main topological ingredient is provided  by the Seiberg-Witten theory of 
the link $M$. 

For any fixed $spin^c$-structure $[k]\in Spin^c(M)$, one defines the 
{\em modified Seiberg-Witten invariant} $\ssw^0_M([k])$ as the sum of the
number of Seiberg-Witten monopoles and the Kreck-Stolz invariant,
see \cite{Chen1,Lim,MW} 
(for this notation, more discussions and references, see \cite{SWI}). 
In this article we prefer to change its sign: we will write 
$\ssw_{M,[k]}^0:=-\ssw_M^0([k])$. 
In general it is very difficult to compute $\ssw_{M,[k]}^0$ 
using its analytic  definition,  therefore  there is an  intense activity
to replace this definition with a different one. 
Presently, there exist a few candidates. In this note we will
discuss two of them. One of them is   $\ssw^{TCW}_{M,[k]}$ provided by 
the sign refined Reidemeister-Turaev torsion (cf. \cite{Tu5})  normalized by the 
Casson-Walker invariant, cf. \ref{5.2}; another is $\ssw^{OSz}_{M,[k]}$
provided by the Ozsv\'ath-Szab\'o theory \cite{OSz}, 
see also \cite{OSzP,NOSZ} for possible connections 
with singularities; cf. with \ref{5.3} here.

In particular, we formulate our conjecture for the  ``symbol'' $\ssw^*_{M,[k]}$,
which can be replaced by any of the above invariants. The conjectured
identity  $\ssw^0_{M,[k]}=\ssw^{TCW}_{M,[k]}=\ssw^{OSz}_{M,[k]}$
makes all the conjectures (for different $*$)  equivalent. Nevertheless,
different realizations might illuminate essentially different aspects. 

\subsubsection{}\label{5.2} $\ssw_{M,[k]}^{TCW}$.
For every $spin^c$-structure $[k]$, Turaev  defines
the {\em sign refined Reidemeister-Turaev torsion }
$\et_{M,[k]}=\sum_{h\in H}\et_{M,[k]}(h) h\in \Q[H]$
(determined by the Euler structure) associated with $[k]$ \cite{Tu5}.
Set $\lambda(M)$
for the {\em Casson-Walker invariant of $M$} (normalized as in \cite{Lescop},
4.7).  Then, one defines 
$$\ssw^{TCW}_{M,[k]}:=-\et_{M,[k]}(1)+\frac{\lambda(M)}{|H|},$$
where 1 denotes the neutral element of the group $H$
(with the multiplicative notation). 

\subsubsection{}\label{5.3} $\ssw_{M,[k]}^{OSz}$.
For any oriented rational homology
3-sphere   $M$ and $[k]\in Spin^c(M)$,
the {\em Ozsv\'ath-Szab\'o $\Z[U]$-module $HF^+(M,[k])$}
(or, the absolutely graded Floer homology) 
was introduced in \cite{OSz} (cf. also with the long list of recent preprints 
of Ozsv\'ath and Szab\'o). This has a $\Q$-grading compatible with 
the $\Z[U]$-action, where $\deg(U)=-2$.  

Additionally, $HF^+(M,[k])$ has an (absolute) $\Z_2$-grading; 
$HF^+_{even}(M,[k])$, respectively  
$HF^+_{odd}(M,[k])$ denote the part of 
$HF^+(M,[k])$ with the corresponding parity.   

For any rational number $r\in \Q$,
we denote by 
 $\calt_r^+$ (cf. \cite{OSzP,NOSZ}) the 
graded $\Z[U]$-module defined as the quotient of $\Z[U,U^{-1}]$
by the submodule  $U\cdot \Z[U]$. 
This has a rational grading in such a way that $\deg(U^{-d})=2d+r$ ($d\geq 0$).

By the general theory, one has a graded $\Z[U]$-module isomorphism 
$$HF^+(M,[k])=\calt^+_{d(M,[k])}\oplus HF^+_{red}(M,[k]),$$
(for some well-defined rational number 
$d(M,[k])$), where $HF^+_{red}$ has a finite $\Z$-rank, and an induced 
(absolute) $\Z_2$-grading.
From this one extracts two numerical invariants: $d(M,[k])$ and 
$$\chi(HF^+(M,[k])):=\rank_\Z HF^+_{red,even}(M,[k])
-\rank_\Z HF^+_{red,odd}(M,[k]).$$
%
%
Then one defines 
$$\ssw^{OSz}_{M,[k]}:=\chi(HF^+(M,[k]))-\frac{d(M,[k])}{2}.$$

\subsection{Conjecture (strong version).}\label{5.4} {\em Let $\x$ be a normal surface
singularity whose link $M$ is a rational homology sphere.
Let $\pi:\tx\to X$  be  a fixed good resolution and $s:=\#\calj$ 
the number of irreducible exceptional divisors of $\pi$. 
Consider an arbitrary $l'\in \bL'$ and define the 
characteristic element $k:=K-2l'\in Char$. Then 

(a) For any line bundle $\call\in Pic(\tx)$ with $c_1(\call)=l'$ one has 
\begin{equation*}
h^1(\call)\leq -\ssw^*_{M,[k]}-\frac{k^2+s}{8}.\tag{1}
\end{equation*}

(b) Additionally, if $\call=\calo_{\tx}(l')$ and 
$\x$ is $\Q$-Gorenstein (i.e. some power of the 
line bundle of holomorphic 2-forms over $X\setminus \{0\}$ is
analytically trivial), then in (1) one has equality.}

(For a generalization see \ref{5.5}(d),
 where the restriction $l'\in \bL'$ is dropped.)

\subsection{Remarks.}\label{5.5}  \

(a) If $\call=\calo_{\tx}$ then $h^1(\calo_{\tx})$ is the 
geometric genus $p_g\x$ of $\x$, hence this case corresponds exactly
to the conjecture formulated in \cite{SWI}.    For detailed historical
 remarks and list of cases   for which the conjecture (for $\calo_{\tx}$)
was verified, see \cite{SWI,SWII,SWIII,NOSZ}. 

(b) Notice that $\ssw^*_{M,[k]}$ depends only on the class $[k]$ of $k$.
In particular, the right hand side of (1) consists of the ``periodical'' term 
$-\ssw^*_{M,[k]}$ and the  ``quadratic'' term  $-(k^2+s)/8$.

(c) In order to prove the conjecture, it is enough to verify it for line bundles $\call$
with $c_1(\call)=l'$ of type  $l'=l'_{ne}(h)$ (for some $h\in H$).

 Indeed, write $l'$ in the form $l'=l'_1+l$ where $l'_1=e(l')=l'_{ne}(l'+L)$ and $l\in L_e$.
Let $RHS(l')$, resp. $RHS(l'_1)$, be the right hand side 
of (1) for $l'$, resp. $l'_1$. Since $[K-2l']=[K-2l'_1]$, the Seiberg-Witten invariants 
are the same, hence
$$RHS(l')-RHS(l'_1)=\frac{-(K-2l')^2+(K-2l'_1)^2}{8}=-(l,l')-\chi(l).$$
This combined with \ref{4.4}(c) shows that
the statements of \ref{5.4} for $\call$  and $\call\otimes
\calo_{\tx}(-l) $ are equivalent. 

(c') Consider any set of representatives $\{l'\}_{l'\in R}$ $(R\subset \bL')$
of the classes $H$, i.e. $\{l'+L\}_{l'\in R}=H$. 
Then the comparison (c) can be also be done for any $\call$ with $l':=
c_1(\call)\in R$ and for $\call\otimes\calo_{\tx}(-l)$.
Therefore, 
the validity of the conjecture \ref{5.4}  follows from the verification of
(1) for line bundles $\call$ with $c_1(\call)\in R$.  

E.g., one can take $R=Q$, or $R=-Q$.
The importance of $R=-Q$ is emphasized by \ref{3.9}.
This fact is exploited in second version of the conjecture.

(d) Similarly, if one verifies the inequality (1) for any $l'\in \bL'$,
then one gets automatically the inequality (1) for {\em any}
$l'\in L'$, hence for any $\call\in Pic(\tx)$.
This statement follows by induction:
if the inequaity (1) is valid for some $\call$, then it is valid for 
$\call\otimes \calo_{\tx}(-E_j)$ (for any $j\in \calj$). 

Indeed, using the exact sequence $0\to 
\call\otimes \calo_{\tx}(-E_j)
\to \call \to \call|_{E_j}\to 0$, one gets that
$h^1(\call\otimes \calo_{\tx}(-E_j))\leq h^1(\call)+1+(c_1(\call),E_j)$.
The proof ends with similar comparison as in (c). 

The point is that if $l'\not\in \bL'$, then the inequality (1),
in general, is not sharp (optimal). 

\subsection{Example. Almost-rational singularities.}\label{alrat} \cite{NOSZ}
First we recall that a singularity is rational if its graph is rational, 
namely, if  for any $l>0$ one has $\chi(l)>0$ (cf. \cite{Artin62,Artin66}).
We say that $\x$ is {\em almost-rational} 
(in short {\em AR}) if one of its good resolution graphs is almost-rational.
A graph $\Gamma$ is almost rational if it is a negative definite 
connected tree which has a distinguished 
vertex $j_0\in\calj$ such that replacing $e_{j_0}$ by some  
$e_{j_0}'\leq e_{j_0}$  we get a rational graph.

The set of {\em AR} graphs is rather large (cf. \cite{NOSZ}
8.2): all rational, weakly elliptic and star-shaped
 graphs are {\em AR}. Moreover, they contain all 
the graphs considered  by Ozsv\'ath and Szab\'o in \cite{OSzP}.
The class of {\em AR} graphs is closed while taking subgraphs and 
decreasing the Euler numbers $e_j$.

The main point of \cite{NOSZ} is that for 3-manifolds with {\em AR} 
plumbing graphs, the Seiberg-Witten invariants can be computed 
combinatorially using ``computational sequences'' (objects in general used in 
singularity theory to evaluate $p_g$). In fact, theorem 9.6 of \cite{NOSZ}
implies the following.

\subsection{Fact.}\label{Fact} {\em Part 
(a) of the strong version of the conjecture (hence \ref{5.5}(d) as well)
 is true for any almost rational singularity.} 

\subsection{Discussion.}\label{dis} (a) Let $(X_a,0)$ be the universal abelian cover
of $\x$ with its natural $H$-action. Obviously, if $p:\tz\to X_a$ is a 
resolution of $(X_a,0)$, then $\tz$ inherits a natural $H$-action.
Recall that the geometric genus $p_g(X_a,0)$ of $(X_a,0)$ is 
$h^1(\tz,\calo_{\tz})$. But one can define much finer invariants:
consider the eigenspace decomposition $\oplus _{\chi\in \hat{H}}H^1
(\tz,\calo_{\tz})_\chi$ of $H^1(\tz,\calo_{\tz})$, and take
$$p_g(X_a,0)_\chi:=\dim_\C H^1(\tz,\calo_{\tz})_\chi \ \  \ 
\mbox{(for any $\chi\in \hat{H}$)}.$$

(b) Obviously, we can repeat the above definition for {\em any } 
(unbranched) abelian cover of
$\x$. More precisely, for any epimorphism $H\to K$ one can
take the composed map $\pi_1(M)\to H\to K$ which defines 
a Galois $K$-covering $(X_K,0)\to \x$ of $\x$ (with $(X_K,0)$ normal). 
Similarly as above,  one can define
$p_g(X_K,0)_\chi$ for any $\chi\in\hat{K}$. But these invariants are not 
essentially new: all of 
them can be recovered  from the corresponding invariants associated with the 
universal abelian cover. Indeed, consider $\chi\in \hat{H}$ via $\hat{K}
\hookrightarrow \hat{H}$. Then $p_g(X_K,0)_\chi=p_g(X_a,0)_\chi$. 
In particular, 
$$p_g(X_K,0)=\sum _{\chi\in\hat{K}}\, p_g(X_a,0)_\chi.$$

(c) In the above definition (part (a)), 
$p_g(X_a,0)_\chi$ is independent of the choice of $\tz$, in particular
one can take $\tz$ considered in the proof of \ref{prop3}.  Those facts,
together with \ref{3.9} show that
$$p_g(X_a,0)_{\theta(h)}=h^1(\calo_{\tx}(-l'_e(h))) \ \ \ \mbox{(for any $h\in H$)}.$$
Since the set $\{-l'_e(h)\}_{h\in H}$ is a set of representatives (cf. \ref{5.5}(c')),
the previous  version \ref{5.4} of the conjecture, {\em restricted to the set
of line bundles of type $\calo_{\tx}(l')$}
is completely equivalent with the following.

\subsection{Conjecture (weak version).}\label{conj2} {\em Let
 $\x$ be a normal surface
singularity whose link $M$ is a rational homology sphere.
Let $\pi:\tx\to X$  be  a fixed good resolution and $s:=\#\calj$ 
the number of irreducible exceptional divisors of $\pi$. 
For any $h\in H=H_1(M,\Z)$ consider the character $\chi:=\theta(h)$
and the characteristic element $k:=K+2l'_e(h)\in Char$. Then  for any $h\in H$
\begin{equation*}
p_g(X_a,0)_{\theta(h)}\leq -\ssw^*_{M,[k]}-\frac{k^2+s}{8}.\tag{2}
\end{equation*}
Additionally, if $\x$ is $\Q$-Gorenstein then in (2) one has equality.}

\subsection{Assumption.}\label{A} For some  $\ssw^*$-theory  and link 
$M$ consider the identity
\begin{equation*}
\sum_{[k]\in Spin^c(M)}\ \ssw^*_{M,[k]}=\lambda(M).\tag{3}\end{equation*}
Notice that for $\ssw^{TCW}$ (and any $M$), (3) is true
because $\sum_{[k]}\et_{M,[k]}(1)=0$. In particular, we expect that 
(3) is valid in any situation. [Obviously, if for some $M$ one verifies
that $\ssw^{OSz}=\ssw^{TCW}$, then (3) will be valid for $\ssw^{OSz}$ as well.
This is the case of Seifert 3-manifolds, see \cite{NOSZ}.] 

\subsection{Corollary.}\label{corab} {\em 

(I). (2) and (3) imply
\begin{equation*}
p_g(X_a,0)\leq -\lambda(M)-\sum_{l'\in Q}\, \frac{(K+2l')^2+s}{8}
\tag{4}\end{equation*}
for any normal singularity $\x$, with equality if $\x$ is $\Q$-Gorenstein.

(II). Assume that for some link $M$ the identity (3) works,
and also the inequality (2) was verified (see e.g. \ref{alrat}).
Then the equality (4) implies the Conjecture \ref{conj2} with equalities 
for all $h$.}

\subsection{Example.}\label{wh} For example, for singularities with good
$\C^*$ action the link is a Seifert manifold. In this case
the inequality in (2) and the identity (3) was verified in \cite{NOSZ}
(cf. also \ref{alrat}). Therefore, the weak version \ref{conj2}
of the conjecture  with equalities follows from  (4) with equality.
We exemplify this by  a Brieskorn singularity.
The general case will be treated in the forthcoming paper \cite{NH}.

Assume that $\x=\{x^2+y^3+z^{12t+2}=0\}$. Here we assume that $t\geq 1$
(if $t=0$ then $\x$ is rational, a case which will be clarified in the next 
section.) The following invariants of $\x$   can be computed using
\cite{SWI}, section 6.

The minimal good resolution graph of $M$ is star-shaped with
three arms and (normalized)  Seifer invariants $(\alpha,\omega)$ equal
to $(3,1), (3,1), (6t+1,2t)$, the self intersection of the central curve
is $-1$, the 
orbifold euler number equals $-1/(18t+3)$, $K^2+s=2$, $\lambda(M)=
-(24t+1)/12$, and $H=\Z_3$. 

Consider the two arms corresponding to Seifert invariants $(3,1)$. Both of 
them contain only one vertex with self intersection $-3$. We denote them 
by $E_1$, respectively $E_2$. Then $Q$ contains exactly three element.
They are $0$, $l_1':=(E_1+2E_2)/3$ and $l_2':=(2E_1+E_2)/3$. One can rewrite
$$\sum_{l'\in Q}\, \frac{(K+2l')^2+s}{8}=|H|\cdot \frac{K^2+s}{8}-
\sum_{l'\in Q} \, \chi(l').$$
In our case, by an easy verification $\chi(l_1')=\chi(l_2')=1/3$. 
Therefore, the right hand side of (4) is $(24t+1)/12-3\cdot 2/8+2/3=2t$.

On the other hand, the universal abelian cover of $\x$ is isomorphic 
to a Brieskorn singularity of type $(X_a,0)=\{u^3+v^3+w^{6t+1}=0\}$
(with the action $\xi*(u,v,w)=(\xi u,\bar{\xi}v,w)$, $\xi^3=1$). 
And one can verify easily that $p_g(X_a,0)=2t$. 
Hence, for $\x$, (2) is valid with equalities.

In fact, in this case $p_g(X,0)=2t$ as well, hence $p_g(X_a,0)_\chi=0$
for any $\chi\not=1$. 

\subsection{}\label{comp} It is also instructive to compare the 
statements of the conjecture
applied for both $\x$ and for some abelian cover $(X_K,0)$ of
$\x$ (cf. \ref{dis}(b)), provided that the link of
$(X_K,0)$ is again a rational homology sphere. 
In the next corollary we exemplify this in a simple
situation, the interested reader is invited to work out more complicated cases.
Write $M_a, \ K_a$ and $s_a$ for the 
corresponding invariants  associated with $(X_a,0)$.
(In particular, $M_a$ is the universal abelian regular cover of $M$.)
 Then one has:

\subsection{Corollary.}\label{cor2} {\em Assume that $\x$ is $\Q$-Gorenstein,
$M$ is a rational homology sphere and $M_a$ is an integral homology sphere.
Then expressing $-p_g(X_a,0)$ in two different ways, one gets:} 
$$\lambda(M_a)+\frac{K_a^2+s_a}{8}=\lambda(M)+
\sum_{l'\in Q}\, \frac{(K+2l')^2+s}{8}.$$

\section{The case of rational singularities}

\subsection{}\label{6.1} In this section we assume that $\x$ is rational, and
$\ssw^*=\ssw^{OSz}$. Some of the arguments are based
on the results of \cite{NOSZ} about $\ssw^{OSz}_{M,[k]}$.
E.g., by \cite{NOSZ}  6.3 and 8.3, one has:
\begin{equation*}
-\ssw^{OSz}_{M,[k]}=\frac{d(M,[k])}{2}=\frac{k_r^2+s}{8}
\ \ \ (\mbox{cf. also with \ref{2.4}}).\tag{1}
\end{equation*}
Recall also that any rational singularity is automatically $\Q$-Gorenstein.

 \subsection{Theorem.}\label{6.2} {\em Conjecture \ref{5.4} 
(hence \ref{conj2} too)  is true for any 
rational singularity if  we take $\ssw^*_{M,[k]}=\ssw^{OSz}_{M,[k]}$.}

\begin{proof}
By \ref{5.5}(c), we can assume  that $l'=l'_{ne}(h) $ where $h=l'+L$.
Using \ref{4.3}, and the fact that $l'_{ne}(h)\in NE_\Q$,
one gets that $h^1(\calo_{\tx}(l'))=0$. On the other hand, by (1),
$-\ssw^{OSz}_{M,[k]}=(k_r^2+s)/8$. 
Here (see also \ref{2.4}), $k=K-2l'$, hence  $k_r=K+2\bar{l}'_{ne}(-l'+L)$ with
$\bar{l}'_{ne}(-l'+L)=-l'_{ne}(-(-l'+L))=-l'_{ne}(l'+L)=-l'$. Hence $k_r=k$ and the right hand
side of \ref{5.4}(1) is also vanishing. 
\end{proof}

Theorem \ref{6.2} (and its proof) and the identity \ref{6.1}(1) 
have  the following consequence (whose
statement is independent of the  Seiberg-Witten  theory):

\subsection{Corollary.}\label{6.3} {\em Assume that $\x$ is rational, and fix 
some $h\in H$. Then
$$p_g(X_a,0)_{\theta(h)}=\frac{(K+2\bar{l}'_{ne}(h))^2-(K+2l'_e(h))^2}{8}=
-\chi(\bar{l}'_{ne}(h))+\chi(l'_e(h)).$$
In particular, $(X_a,0)$ is rational if and only if 
$\chi(\bar{l}'_{ne}(h))=\chi(l'_e(h))$ for any $h\in H$.}

\vspace{2mm}

This emphasizes in an impressive way the differences between the two
``liftings'' $l'_e(h)$ and $\bar{l}'_{ne}(h)$. Recall: for a class $h=l'+L$,
both $l'_e(h)$ and $\bar{l}'_{ne}(h)$ are elements of $l'+L$, but
the first is minimal in $L_{\Q,e}$ (i.e. it is the representative in $Q$),
while the second is minimal in $-NE_\Q$. Since $-NE_\Q\subset 
L_{\Q,e}$, one has $l'_e(h)\leq \bar{l}'_{ne}(h)$. In  some cases they are not equal. 

For example, take the $A_4$ singularity, where $E$ has three components
$E_1,\, E_2, \, E_3$, all  with self-intersection $-2$,  and $E_2$ intersecting
the others transversely. Then $-D_2=(1/2,1,1/2)=\bar{l}'_{ne}(h)$ for some 
$h$, but it is not in $Q$: $l'_e(h)=-D_2-E_2=(1/2,0,1/2)$. 
Nevertheless, in this case, their Euler characteristics are the same
(corresponding to the fact that -- $A_4$ being a cyclic quotient singularity -- the
universal abelian cover is smooth).  

\subsection{Example.}\label{6.4} Assume that $\x$ is a rational singularity with 
the following dual  resolution graph.

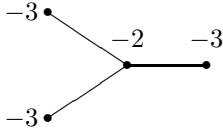
\begin{picture}(300,60)(0,10)
\linethickness{.5pt}
\put(70,40){\circle*{3}}
\put(100,40){\circle*{3}}
\put(40,60){\circle*{3}}
\put(40,20){\circle*{3}}
\put(40,60){\line(3,-2){30}}
\put(40,20){\line(3,2){30}}
\put(70,40){\line(1,0){30}}
\put(70,50){\makebox(0,0){$-2$}}
\put(100,50){\makebox(0,0){$-3$}}
\put(30,60){\makebox(0,0){$-3$}}
\put(30,20){\makebox(0,0){$-3$}}
\end{picture}

Let $\tx$ be its minimal resolution, we write $E_0$ for the central
exceptional component, the others are denoted by $E_i$, $i=1,2,3$.
Then by a computation $l':= -(D_1+D_2+D_3)=(1,2/3,2/3,2/3)$. This element 
is minimal  in $-NE_\Q$ but not in $L_{\Q,e}$. 
Its representative in $Q$  is $l'-E_0$.
The character $\chi=\theta([l'])$ is $\chi(g_0)=1$ and $\chi(g_j)=e^{4\pi i/3}$
for $j=1,2,3$. By a computation $K=-l'$, $(K+2l')^2=-2$ and 
$(K+2(l'-E_0))^2=-10$. Therefore, $p_g(X_a,0)_\chi=(-2+10)/8=1$.
In particular, $(X_a,0)$ cannot be rational.

Indeed, one can verify (using e.g.  \cite{Neu}), that the minimal resolution 
of $(X_a,0)$  contains exactly  one irreducible exceptional curve of genus 1 and 
self-intersection $-3$. In particular, $(X_a,0)$ is minimally elliptic
with $p_g(X_a,0)=1$. This also shows that the above eigenspace is the only non-trivial 
one. 

We reverify this last fact  for the conjugate $\bar{\chi}$ of $\chi$.
In this case $h=\theta^{-1}(\bar{\chi})$ is the class of 
$-D_0=(1,1/3,1/3,1/3)=\bar{l}'_{ne}(h)$; and 
$l'_e(h)=-D_0-E_0$.  By a calculation $K+2\bar{l}'_{ne}(h)=E_0$ and 
$K+2l'_{e}(h)=-E_0$, hence their squares are the same. In particular,
$p_g(X_a,0)_{\bar{\chi}}=0$.


\begin{thebibliography}{99}

\bibitem{Artin62} Artin, M.: 
Some numerical criteria for contractibility of curves on algebraic surfaces.
{\em  Amer. J. of Math.}, {\bf 84} (1962),  485-496.

\bibitem{Artin66} Artin, M.: 
On isolated rational singularities of surfaces.
{\em Amer. J. of Math.}, {\bf 88} (1966), 129-136.

\bibitem{Chen1} Chen, W.: Casson invariant and Seiberg-Witten gauge
theory, {\em Turkish J. Math.}, {\bf 21}(1997), 61-81.

\bibitem{GS} Gompf, R.E. and Stipsicz, I.A.:  An Introduction to
$4$-Manifolds and Kirby Calculus, {\em Graduate Studies in Mathematics}, vol. {\bf 20},
Amer. Math. Soc., 1999.

\bibitem{GRa} Grauert, H.: \"Uber Modifikationen  und exceptinelle
analytische Mengen, {\em Math. Annalen}, {\bf 146} (1962), 331-368.

\bibitem{GR} Grauert, H. and Remmert, R.: Komplexe R\"aume, {\em Math. Ann.}
{\bf 136} (1958), 245-318.

\bibitem{Kollar} Koll\'ar, J,: {\em Shafarevich  Maps and Automorphic Forms},
Princeton University Press, 1995. 

\bibitem{Laufer72} Laufer, H.B.: On rational singularities,
{\em Amer. J. of Math.}, {\bf 94} (1972), 597-608.


\bibitem{Lescop} Lescop, C.: Global Surgery Formula for the Casson-Walker
Invariant, {\em Annals of Math. Studies}, vol. {\bf 140}, Princeton University
Press, 1996.

\bibitem{Lim} Lim, Y.: Seiberg-Witten invariants for 3-manifolds in the
case $b_1=0$ or $1$, {\em Pacific J. of Math.}, {\bf 195}(2000), 179-204.

\bibitem{MW}  Marcolli, M. and Wang, B.L.: Seiberg-Witten invariant and the 
Casson-Walker invariant for rational homology 3-spheres,
math.DG/0101127, Geometri{\ae} Dedicata, to appear.

\bibitem{Ninv}  N\'emethi, A.:
``Weakly'' Elliptic Gorenstein Singularities of
Surfaces, {\em Inventiones math.}, {\bf 137} (1999), 145-167.


\bibitem{NOSZ} N\'emethi, A.: On the Ozsv\'ath-Szab\'o invariant of negative 
definite plumbed 3-manifolds, preprint. 

\bibitem{NH} N\'emethi, A.: Line bundles associated with normal singularities,
II (singularities with $\C^*$ action), in preparation. 

\bibitem{SWI} N\'emethi, A. and Nicolaescu, L.I.: 
 Seiberg-Witten invariants and surface singularities,
{\em Geometry and Topology},  Volume {\bf 6} (2002), 269-328. 

\bibitem{SWII} N\'emethi, A. and Nicolaescu, L.I.: 
 Seiberg-Witten invariants and surface singularities II
(singularities with good $\C^*$-action), 
math.AG/0201120.

\bibitem{SWIII} N\'emethi, A. and Nicolaescu, L.I.: 
 Seiberg-Witten invariants and surface singularities III
(splicings and cyclic covers)
math.AG/0207018.

\bibitem{Neu} Neumann, W.:  Abelian covers of quasihomogeneous surface
singularities, {\em Proc. of Symposia in
Pure Mathematics},  {\bf  40} (2) (1983), 233-244.

\bibitem{NW} Neumann, W. and Wahl, J.: Casson invariant of links
of singularities, {\em Comment. Math. Helv.} {\bf 65}, 58-78, 1991.

\bibitem{NWuj} Neumann, W. and Wahl, J.: 
Complex surface singularities with integral homology sphere links,
math.AG/0301165.



\bibitem{OSz} Ozsv\'ath, P.S. and Szab\'o, Z.: Holomorphic discs
and topological invariants for rational homology three-spheres,
math.SG/0101206.




\bibitem{OSzP} Ozsv\'ath, P.S. and Szab\'o, Z.: On the Floer 
homology of plumbed three-manifolds, math.SG/0203265.

\bibitem{Pinkh} Pinkham, H.: {\sl Normal surface singularities with ${\C}^*$
 action}, Math. Ann. {\bf 117}(1977), 183-193.

\bibitem{Reid} Reid, M.: Chapters on Algebraic Surfaces, in: Complex
Algebraic Geometry, IAS/Park City Mathematical Series,
Vol. 3. (J. Koll\'ar editor) 1997.

\bibitem{Sakai} Sakai, F.: Anticanonical models of rational surfaces,
{\em Math. Ann.}, {\bf 269} (1984), 389-410.


\bibitem{Tu5} Turaev, V.G.:   Torsion invariants of $Spin^c$-structures 
on $3$-manifolds, {\em Math. Res. Letters}, {\bf 4}(1997), 679-695.

\end{thebibliography}
\end{document}